\documentclass[conference,12pt,onecolumn]{IEEEtran}

\newcommand{\subparagraph}{}
\usepackage{titlesec}

\usepackage[english]{babel}
\usepackage[utf8x]{inputenc}
\usepackage{bm}
\usepackage{enumitem}
\usepackage{booktabs}
\usepackage[table,xcdraw]{xcolor}

\usepackage{amsmath}
\usepackage{amsfonts}
\usepackage{comment}
\usepackage[lined]{algorithm2e}
\usepackage{amsthm}
\usepackage{graphicx}
\usepackage{lettrine}
\usepackage[colorinlistoftodos]{todonotes}
\usepackage[colorlinks=true, allcolors=black]{hyperref}

\newtheorem{thm}{Theorem}
\newtheorem{rmk}{Remark}
\theoremstyle{definition}

\newcommand*{\QEDB}{\hfill\ensuremath{\square}}%
\newcommand{\bd}{\boldsymbol}

\makeatletter
\renewcommand\paragraph{\@startsection{paragraph}{4}{\z@}%
    {1ex \@plus0ex \@minus0ex}%
    {-1em}%
    {\normalfont\normalsize\bfseries}}
\makeatother

\usepackage{footnote}
\makesavenoteenv{tabular}
\makesavenoteenv{table}
\ifCLASSINFOpdf
 
\fi

\begin{document}

\title{Non-Wire Alternatives to Capacity Expansion}

\author{\IEEEauthorblockN{Jesus E. Contreras-Oca\~na\textsuperscript{a,b}, Uzma Siddiqi\textsuperscript{b}, and Baosen Zhang\textsuperscript{a}} 
\IEEEauthorblockA{ \textsuperscript{a}University of Washington Electrical Engineering, \textsuperscript{b}Seattle City Light \\ 
Emails: \{jcontrer, zhangbao\}@uw.edu, uzma.siddiqi@seattle.gov }}

\maketitle

\begin{abstract}
Distributed energy resources (DERs) can serve as non-wire alternatives to capacity expansion by managing peak load to avoid or defer traditional expansion projects. In this paper, we study a planning problem that co-optimizes DERs investment and operation (e.g., energy efficiency, energy storage, demand response, solar photovoltaic) and the timing of capacity expansion. We formulate the problem as a large scale (in the order of millions of variables because we model operation of DERs over a period of decades) non-convex optimization problem. Despite its non-convexities, we find  its optimal solution by decomposing it using the Dantzig-Wolfe Decomposition Algorithm and solving a series of small linear problems.  Finally, we present a real planning problem at the University of Washington Seattle Campus.  
\end{abstract}

\IEEEpeerreviewmaketitle

\section{Introduction}
Electric utility distribution systems are typically designed and built for peak load which usually happens a small number of hours per year. When the system load reaches capacity, the traditional solution has been to install more wires or reinforce existing ones~\cite{seifi2011electric}. While decades of experience make this solution reliable and safe, it is often associated with enormous capital costs, hostile public opinion, and/or time-consuming legal issues (e.g., eminent domain questions)~\cite{stanton2015getting}. 

Lately, there has been an increased interest in distributed energy resources (DERs) such as energy storage (ES), energy efficiency (EE), demand response (DR), and distributed generation (DG) as alternatives to traditional ``wire" solutions. In the planning community, these solutions are often called non-wire alternatives (NWAs)\footnote{We use the term NWA and DER interchangeably. We employ the term NWA to emphasize their impact on traditional capacity expansion solutions.}. The basic premise is that NWAs can manage load to avoid, or at least delay, the need for traditional capacity expansion. 

\begin{figure} 
\centering
\includegraphics[width=.75\textwidth]{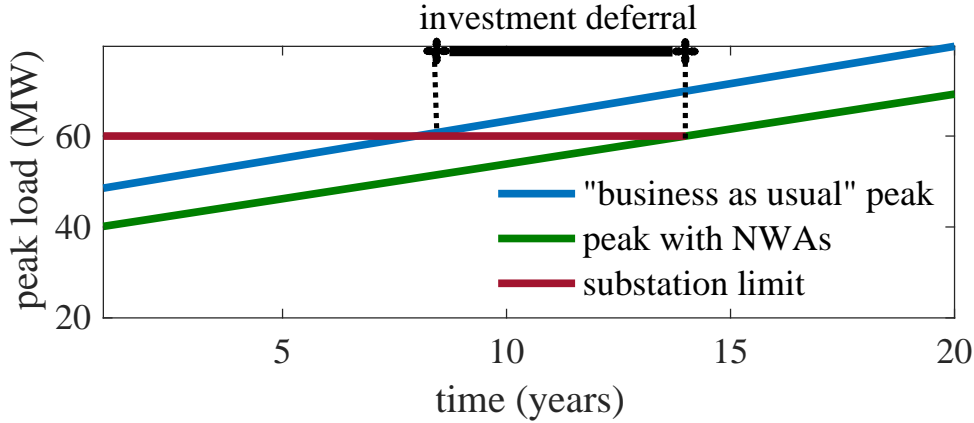}
\caption{Projected load growth at the University of Washington. When the load reaches the limit, the substation capacity must be expanded. Note that capacity expansion can be deferred by managing load growth via NWAs. } \label{fig:upgrade_problem_rev}
\end{figure}

The University of Washington (UW) is expected to add 6 million sq. feet of new buildings (e.g., labs, classrooms, office space) to its Seattle Campus during the next 10 years~\cite{UWMP_2018}. This would translate into approximately 17 MW of additional load and require the capacity at the substation that serves the campus to be expanded. The blue line in Fig.~\ref{fig:upgrade_problem_rev} is the projected ``business as usual'' campus peak load while the green line is peak load managed via a set of NWAs. When the load reaches the feeder or substation limit, the system planner must expand its capacity. In this case, NWAs are able to \emph{delay} the need for capacity expansion by reducing peak load. 

The economic reason for deferring investments is the time-value of money, which states that a dollar spent now is more valuable than a dollar spent later. Policy-wise, there are often other benefits of deferring capital-intensive projects (e.g., reducing the risk of expected load not materializing, generating local employment opportunities, among others~\cite{stanton2015getting}). In this paper, we focus on the economic question and ask: is deferring traditional expansion investments worth the costs of NWAs?

The answer to this question is non-trivial. For one, the cost and benefits of NWAs are not only a function of their installed capacities but also of their operation. Thus, one must co-optimize investment and operation of NWAs to find optimal decisions. This leads to a large problem that can be hard to solve. Furthermore, the considering the time-value of deferring investments introduces non-linearities and non-convexities to the planning problem, making it even harder to solve.
 \subsection{Contributions}
In this paper, we make the following contributions:
\begin{enumerate}
\item A \textbf{formulation of the NWAs planning problem} that determines 1) investment and 2) operation of NWAs and 3) the timing of capacity expansion.
\item A \textbf{scalable solving technique}. The NWAs planning problem is a large scale, non-convex optimization problem. It is large-scale (on the order of millions of variables) because we model the operation of the NWAs over an investment horizon of decades. The problem's non-convexities are introduced by the variables and constraints  used to model the timing of capacity expansion. We deal with the scale of the problem by decomposing it into smaller subproblems using the Dantzig-Wolfe Decomposition Algorithm (DWDA). As shown shortly, the non-convexities end up being confined to the DWDA's master problem which, in our case is small (in the order of tens to a couple of hundred variables). We deal with the non-convexities of the master problem by further decomposing it and solving a small number of linear programs. 
\item A \textbf{case study} where NWAs may be used to defer substation and feeder upgrades at the UW Seattle Campus.  
\end{enumerate}
\subsection{Literature review}
The idea of deferring infrastructure investments by reducing load was first introduced in~\cite{Xin_One_2005}. The work in~\cite{Xin_One_2005} focuses on quantifying the effects of load reduction on avoided infrastructure costs and not on finding the optimal load reduction or the appropriate technologies to do so. On a similar note, the authors of~\cite{Gin_Quant_2006} and~\cite{Piccolo_Evaluating_2009} develop frameworks to quantify the value of capacity deferral of distributed generation by explicitly modeling DG as the mechanism of net load reduction. However, they also do not address the problem of finding optimal DG investment nor consider other types of DERs. In~\cite{Samper_Investment_2013_I}, the authors determine optimal investments in DG considering the value of network investment deferral. However, their non-linear mixed-integer formulation is intractable in general. In contrast, our model considers a wider set of NWAs and is tacked by solving a series of smaller convex optimization problems. 

Furthermore, there is relatively little literature on holistic DER planning. Most consider a narrow definition of the term DER that only includes DG, e.g.,~\cite{Alarcon-Rodriguez_MO_2009, Yuan_CO_2017}, or only ES and DR~\cite{Dvorkin_Merchant_2017}. Instead, we consider a generic definition of DERs and present a case that considers solar photovoltaic (PV) generation, DR, EE, and ES. 

\subsection{Organization of this paper}

This paper is organized as follows. Section~\ref{sec:Substation_upgrade_problem} introduces the capacity expansion problem and formulates it as an optimization problem Section~\ref{sec:NWA} introduces a generic NWA model and four specific instances: EE, PV, DR, and ES. Section~\ref{sec:NWA_PP} introduces the NWAs planning problem and Section~\ref{sec:DWDA} proposes a solution technique. Section~\ref{sec:case_study} presents a case study of load-growth at the University of Washington. Section~\ref{sec:conclusion} concludes the paper.

\section{The Capacity Expansion Problem}
\label{sec:Substation_upgrade_problem}
System planners typically like to expand capacity at the \emph{latest} possible time but in time to meet expected growth. A reason for this is the time-value of money: we would like to spend a dollar later rather than now. Let $l_a^\mathrm{p}$ denote expected peak load\footnote{The expected peak load $\boldsymbol{l}^\mathrm{p}$ is usually calculated by the utility using population growth projections, planned construction projects, weather forecasts, among others~\cite{seifi2011electric}.} during year $a$ and the pre-expansion capacity as $\overline{l}$. After expansion, we assume that any reasonable load can be accommodated for the foreseeable future.  Then, the decision rule for choosing a year to expand capacity is 
\begin{equation} \label{eq:subs_upgrade_prob}
\mathrm{CapEx}(\boldsymbol{l}^\mathrm{p})=a\; | \;l^\mathrm{p}_{a+1}>\overline{l}, \;l_k^\mathrm{p} \le \overline{l}  \;\forall \; k \le a
\end{equation}
where $\boldsymbol{l}^\mathrm{p}=\{l^\mathrm{p}_a\}_{a\in\mathcal{A}}$ and $\mathcal{A}$ is the planning horizon. The decision rule $\mathrm{CapEx}$ states that the planner expands capacity at a future year $a$ immediately before the limit $\overline{l}$ is \emph{first} reached by the load (as illustrated in Fig.~\ref{fig:upgrade_problem_rev}).  In this paper, we analyze capacity expansion at a single point in a radial system (e.g., a feeder or substation) and assume that the downstream network is non-congested. This assumption allows us to disregard the network.

Let  $I$ denote the inflation-adjusted cost of capacity expansion\footnote{We assume that the inflation-adjusted cost of capacity expansion is constant throughout the planning horizon.}. Then, if system capacity is expanded at year $\mathrm{CapEx}(\boldsymbol{l}^\mathrm{p})$, the present cost of the investment is
\begin{equation}
\tilde{I}(\boldsymbol{l}^\mathrm{p}) = \frac{I} {\left(1+\rho\right)^{\mathrm{CapEx}(\boldsymbol{l}^\mathrm{p})}}  \label{eq:IW_PV}
\end{equation}
where $\rho$ is the annual discount rate, a quantity closely related to the real interest rate~\cite{seifi2011electric}.

\begin{rmk} \label{thm:IW_PV_CVX}
Let the series of yearly peak load during the planning horizon $\bd {\mathcal A}$ be denoted as $\bd l^\mathrm{p} =\{l_a^\mathrm{p}\}_{a \in \mathcal A}$. The function $\tilde{I}$ from Eq.~\eqref{eq:IW_PV} can be reformulated as the optimization problem
\begin{subequations}  \label{eq:IW_PV_CVX}
	\begin{align}
\tilde{I}(\bd{l}^\mathrm{p})=\min_{\delta}& \frac{I}{(1+\rho)^{\delta}}  \\
\mbox{s.t. } & 0 \le \delta \le |\mathcal A| \\
& l_{a}^\mathrm{p} \le \overline{l}\; \forall \;\;\; a<\delta. \label{eq:NCVX_Const}
\end{align}
\end{subequations}
\end{rmk}
As we show in Theorem~\eqref{eq:CVX_problem}, the reformulation of $\tilde{I}$ allows us to formulate a NWA planning problem which simultaneously optimizes NWA investment, operation, and the timing of capacity expansion) with a structure that is conductive to be decomposed into smaller problems.

Note that given a $\bd{l}^\mathrm{p}$,  Problem~\eqref{eq:IW_PV_CVX} is convex. However, if we treat $\bd{l}^\mathrm{p}$ as a variable, its feasible solution space is non-convex. This fact is relevant on a NWA planning context where we can manipulate the peak load and thus treat $\bd{l}^\mathrm{p}$ as a variable. While it is unfortunate that Problem~\eqref{eq:IW_PV_CVX} is non-convex on $\bd{l}^\mathrm{p}$, Section~\ref{sec:DWDA} shows that we can handle its non-convexities by solving at most $|\mathcal A|$ small-scale linear problems\footnote{The length of the planning horizon, $|\mathcal A|$, is in the order tens of years.}  that can be quickly and reliably solved using off-the-shelf solvers. In the next Section, we introduce models of generic and specific NWAs and discuss how they relate to the capacity expansion problem. 

\section{Non-wire alternatives}
\label{sec:NWA}
 Let the index $i$ denote a NWA technology. A generic NWA model is characterized by six elements:
\begin{enumerate}[noitemsep,topsep=0pt]
\item investment decision variables $\phi_i$,
\item operating decision variables $x_i$, 
\item a set of feasible investment decisions $\bd \Phi_i$, 
\item a set of feasible operating regimes $\bd{\mathcal{X}}_i(\phi_i)$,
\item a set of functions $l^i_{a,t}(x_i)$ that map operating decisions onto load at time $t$ of year $a$, and 
\item an investment cost function $I^\mathrm{NW}_i(\phi_i)$. 
\end{enumerate}

While investment decisions are made in horizons on the order of years, operating decisions are made on much shorter horizons. In this paper, operating decisions are made in $\Delta t$-hour intervals (e.g., 1-hour interval in our case study). The set $\mathcal{T}$ denotes the set of operating time intervals during one year.

 In this paper, the feasible solution spaces of investment and operating decisions, $\bd \Phi_i$ and $\bd{\mathcal{X}}_i(\phi_i)$ respectively, are assumed to be convex on both $\phi_i$ and $x_i$. Additionally, we assume that the investment cost function is convex and that the load functions $l^i_{a,t}(x_i)$ are linear in $x_i$. Now we define each of the six elements that define NWAs for the four particular technologies considered in this paper: EE, PV, DR, and ES.

\paragraph*{Energy Efficiency} We model EE as percentage reduction with respect to a base load that translates into a load reduction of $r_{a,t}^\mathrm{EE}$ MWs at all time periods $t$ of all years $a$. For EE, the investment decision is to choose a load reduction percentage. We model the investment cost, $I^\mathrm{NW}_\mathrm{EE}$, as a convex piece-wise linear function of the load reduction percentage~\cite{brown2008us}. The slope of each of the segment $B^\mathrm{EE}$ segments, $C_b^\mathrm{EE}$, represents the marginal cost of load reduction. The six parameters that define DR as a NWA are
\begin{subequations}
\begin{align*}
&\phi_\mathrm{EE} = \left\{ \epsilon_b^\mathrm{EE}\right\}_{b=1,\hdots,B^\mathrm{EE}} \\
& x_\mathrm{EE} = \{r^\mathrm{EE}_{a,t}\}_{a\in \mathcal{A},\; t \in \mathcal{T}}\\
&\Phi_\mathrm{EE} = \left\{ \epsilon_b^\mathrm{EE} \; | \;  \epsilon_b^\mathrm{EE}\in \left[ 0, \overline{\epsilon}_b^\mathrm{EE}\right] \; \forall\; b=1\hdots, B^\mathrm{EE} \right\} \\
&\mathcal{X}_\mathrm{EE} \left(\phi_\mathrm{EE}\right) = \{r^\mathrm{EE}_{a,t} \; |\; r^\mathrm{EE}_{a,t} =l_{0,t}^\mathrm{b} \cdot \sum_{b=1}^{B^\mathrm{EE}} \epsilon^{EE}_b \;\forall \; a\in \mathcal{A},\; t \in \mathcal{T}\}\\
&l_{a,t}^\mathrm{EE}(x_\mathrm{EE}) = -r_{a,t}^\mathrm{EE} \\
&I^\mathrm{NW}_\mathrm{EE}(\phi_\mathrm{EE}) =   \sum_{b=1}^\mathrm{B^\mathrm{EE}} C^\mathrm{EE}_b \cdot \epsilon^\mathrm{EE}_b 
\end{align*}
\end{subequations}
where $\epsilon_b^\mathrm{EE}$ is the percentage reduction for each cost piece-wise linear segment of $I^\mathrm{NW}_\mathrm{EE}$, $\overline{\epsilon}_b^\mathrm{EE}$ is the size of each segment, and $l_{0,t}^\mathrm{base}$ is the base load (i.e., the pre-EE load).

\paragraph*{Solar photovoltaic generation} 

For the solar PV case, investment decision $\phi_\mathrm{PV}$ is the PV installed capacity $g^\mathrm{PV}_\mathrm{CAP}$.  At a given time $a,t$, the solar energy generation is given by $g^\mathrm{PV}_{a,t} = \alpha_{a,t}^\mathrm{PV}\cdot g^\mathrm{PV}_\mathrm{CAP}$ where $\alpha_{a,t}^\mathrm{PV}$ is a number in $[0,1]$ and is related to solar radiation levels. All in all, the parameters that define solar PV as a NWA are
\begin{subequations}
\begin{align*}
& \phi_\mathrm{PV}  =  g^\mathrm{PV}_\mathrm{CAP} \\
& x_\mathrm{PV} = \{g_{a,t}^\mathrm{PV}\}_{a\in \mathcal{A},\; t \in \mathcal{T}} \\
&\Phi_\mathrm{PV} = \left\{ g^\mathrm{PV}_\mathrm{CAP}\; |\; g^\mathrm{PV}_\mathrm{CAP} \in \left[0, \overline{g}^\mathrm{PV}_\mathrm{CAP} \right] \right\} \\ 
& \mathcal{X}_\mathrm{PV} \left(\phi_\mathrm{PV}\right) = \{g_{a,t}^\mathrm{PV} \; |\;g_{a,t}^\mathrm{PV} = \alpha_{a,t}^\mathrm{PV} g^\mathrm{PV}_\mathrm{CAP}  \;\forall \; a\in \mathcal{A},\; t \in \mathcal{T}\} \\
& l_{a,t}^\mathrm{PV}(x_\mathrm{PV}) = -g_{a,t}^\mathrm{PV}\\
&I^\mathrm{NW}_\mathrm{PV}(\phi_\mathrm{PV}) = C^\mathrm{PV} \cdot g^\mathrm{PV}_\mathrm{CAP}
\end{align*}
\end{subequations}
where $C^\mathrm{PV}$ is the cost per unit capacity of solar PV (capital and labor costs) and $\overline{g}^\mathrm{PV}_\mathrm{CAP}$ is the limit on PV installed capacity. 

\paragraph*{Demand response}
We consider investments in DR communication and control infrastructure that allow a portion of the load, (e.g.,  water heaters and HVAC systems), to be shifted in time. For DR, the investment decision is the amount DR-enabled load, $r^\mathrm{DR}_\mathrm{CAP}$ which limits the demand reduction $r_{a,t}^\mathrm{DR}$ that can be deployed at a particular time. Since DR allows load to be shifted in time a load reduction of $r_{a,t}^\mathrm{DR}$ at a time $a,t$ is associated with a demand rebound of $\alpha^\mathrm{DR}\cdot r_{a,t}^\mathrm{DR}$ during the next time period. The coefficient $\alpha^\mathrm{DR}$ is a number $\ge 1$ and is related to efficiency losses due to DR deployment. Note that more sophisticated rebound models such as the ones in~\cite{lutolfimpact} are admissible in our framework. The parameters that define DR as a NWA are
\begin{subequations}
\begin{align*}
&\phi_\mathrm{DR}  =  r^\mathrm{DR}_\mathrm{CAP} \\
& x_\mathrm{DR} = \{r_{a,t}^\mathrm{DR}\}_{a\in \mathcal{A},\; t \in \mathcal{T}}\\
& \Phi_\mathrm{DR} = \left\{ r^\mathrm{DR}_\mathrm{CAP}\; | \; r^\mathrm{DR}_\mathrm{CAP} \in\left[ 0, \overline{r}^\mathrm{DR}_\mathrm{CAP}  \right] \right\} \\
&\mathcal{X}_\mathrm{DR} \left(\phi_\mathrm{DR}\right) = \{ r_{a,t}^\mathrm{DR} |  r_{a,t}^\mathrm{DR} \in\left[ 0, r^\mathrm{DR}_\mathrm{CAP} \right] \; \forall \; a\in \mathcal{A},\;t \in \mathcal{T} \} \\ 
& l_{a,t}^\mathrm{DR}(x_\mathrm{DR}) =\alpha^\mathrm{DR}r_{a,t-1}^\mathrm{DR} -r_{a,t}^\mathrm{DR} \\ 
&I^\mathrm{NW}_\mathrm{DR}(\phi_\mathrm{DR}) = C^\mathrm{DR} r^\mathrm{DR}_\mathrm{CAP} 
\end{align*}
\end{subequations}
where $C^\mathrm{DR}$ is the cost of enabling DR per unit capacity of load.

\paragraph*{Lithium-ion energy storage} 
The investment decision for ES is the energy storage capacity $s^\mathrm{max}_0$ of the system. The operating variables are the ES's charge $c_{a,t}$, discharge $d_{a,t}$, state-of-charge $s_{a,t}$, and storage capacity\footnote{$s^\mathrm{max}_a$ may be different than $s^\mathrm{max}_0$ because we model battery degradation. Although $s^\mathrm{max}_a$ is not normally considered an operating variable we do so for ease of notation.} during year $a$, $s^\mathrm{max}_a$. The feasible region of the operating variables is described by Equations~\eqref{eq:feas_reg_start}-\eqref{eq:feas_reg_end} and includes the usual charge, discharge, and state-of-charge limits used to model ES~\cite{Sarker_Optimal_2017,Dvorkin_Merchant_2017}. Additionally, as expressed in Equation~\eqref{eq:degradation}, the storage capacity degrades by $\beta^\mathrm{ESD}$ per unit charge/discharge~\cite{Sarker_Optimal_2017}. The parameters that define ES as a NWA are
\begin{subequations}
\begin{align}
&\phi_\mathrm{ES}  =   s^\mathrm{max}_0 ,\; x_\mathrm{ES} = \{c_{a,t},\; d_{a,t},\; s_{a,t} \\
& s_a^\mathrm{max}\}_{a\in \mathcal{A},\; t \in \mathcal{T}}\\
&\Phi_\mathrm{ES} =  \left\{  s^\mathrm{max}_0 \;|\; s^\mathrm{max}_0 \in \left[0,\overline{s}^\mathrm{max}_0 \right] \right\}\\
&\mathcal{X}_\mathrm{ES} \left(\phi_\mathrm{ES}\right) = \biggl\{c_{a,t},\; d_{a,t},\; s_{a,t},\; s_a^\mathrm{max}\;| \label{eq:feas_reg_start} \\
& s_{a,t+1} = s_{a,t} +  \eta_c \Delta t c_{a,t} -\frac{\Delta t d_{a,t}}{\eta_d} \; \forall \; a\in \mathcal{A},\; t \in \mathcal{T} \\ 
 &s_{a}^\mathrm{max} = s^\mathrm{max}  - \beta^\mathrm{ESD} \cdot\sum_{k = 1}^{a-1} \sum_{t \in \mathcal{T}} (c_{k,t} + d_{k,t}) \; \forall \; a\in\mathcal{A}  \label{eq:degradation} \\
 & s_{a,t} \in \left[0, s_a^\mathrm{max}\right],\;c_{a,t} ,d_{a,t} \in \left[0, \frac{s^\mathrm{max}}{\alpha^\mathrm{EPR}}\right] \; \forall \; a\in \mathcal{A},\; t \in \mathcal{T}  \biggr\} \label{eq:feas_reg_end} \\
 &l_{a,t}^\mathrm{ES}(x_\mathrm{ES}) =c_{a,t} -d_{a,t} \\
 &I^\mathrm{NW}_\mathrm{ES}(\phi_\mathrm{ES}) = C^\mathrm{ES} \cdot s^\mathrm{max}_0
\end{align}
\end{subequations}
where $\eta_c \; (\eta_d)$ is the charge (discharge) efficiency, $\alpha^\mathrm{EPR}$ denotes the energy-to-power ratio of the ES system, and dollar per unit energy cost of $C^\mathrm{ES}$ is the dollar per unit energy cost of storage capacity.  In this work, we consider investments in lithium-ion ES although other chemistries are compatible within our framework~\cite{Sarker_Optimal_2017}.

\begin{rmk} Our framework allows for other NWAs to be included, e.g., controllable electric-vehicles, diverse battery chemistries, dispatchable DG, etc. This is especially true since our solving method, discussed in Section~\ref{sec:NWA_PP}, allows for parallel computation of investment and operating decisions of each of the NWAs. 
\end{rmk}

\subsection{The substation upgrade problem revisited}

Let the total load at time $t$ of year $a$ (with a set $\mathcal{N}$ of NWAs) be denoted by $l_{a,t}(\boldsymbol{x}) = l^\mathrm{b}_{a,t} + \sum_{i \in \mathcal{N}} l_{a,t}^i(\boldsymbol{x}_i)$ where $\boldsymbol{x} = \{x_i \}_{i\in \mathcal{N}}$. Then, the yearly peak load as a function of NWA operation is $l_a^\mathrm{p}(\boldsymbol{x}) = \max_{t \in \mathcal{T}}\{ l_{a,t}(\boldsymbol{x}) \}$. 

Recall from Eqs.~\eqref{eq:IW_PV} and~\eqref{eq:subs_upgrade_prob} that the system planner only had to decide on when to invest in capacity expansion. With NWAs, however, the present cost of expansion
\begin{equation*}
\tilde{I}(\boldsymbol{l}^\mathrm{p}(\boldsymbol{x}))
\end{equation*}
becomes a function of the NWAs operating variables, giving the system planner has the ability to plan investment and operation of NWAs that minimize $\tilde{I}(\boldsymbol{l}^\mathrm{p}(\boldsymbol{x}))$. However, a good plan should also consider the investment cost of the NWAs, their operating costs (e.g., energy costs), and benefits other than deferring capacity expansion (e.g., demand charge reductions). In the next Section, we present a holistic NWA investment/operation and timing of capacity expansion problem.

\section{The Non-wire alternatives planning problem}
\label{sec:NWA_PP}
We would like to determine investments and operation of NWAs and the timing of capacity expansion that minimize total present cost composed of: 1) convex operating cost functions for each NWA, $C_i^\mathrm{O}(x_i)$, 2) NWA investment costs, $I_i^\mathrm{NW}(\phi_i)$, a convex  peak demand charge, $C^\mathrm{D}(\boldsymbol{l}^\mathrm{p}(\boldsymbol{x}))$, and the 3) present cost of capacity expansion, $\tilde{I}(\boldsymbol{l}^\mathrm{p}(\boldsymbol{x}))$. We formally state the NWAs planning problem as
\begin{equation}
\min_{\substack{\phi_i \in \Phi_i \\ x_i \in \mathcal{X}_i (\phi_i) }} \biggl\{\sum_{i \in \mathcal{N}}  \left[ C_i^\mathrm{O}(x_i) + I_i^\mathrm{NW}(\phi_i) \right]  +C^\mathrm{D}(\boldsymbol{l}^\mathrm{p}(\boldsymbol{x}))  + \tilde{I} (\boldsymbol{l}^\mathrm{p}(\boldsymbol{x})) \biggr\} .\label{prob:NWA_planning}
 \end{equation}

From the definition of $\tilde I$ in Eq.~\eqref{eq:IW_PV}, contains non-convex function $\mathrm{CapEx}$ that is inconvenient to use in large-scale optimization problems. In the following theorem, we demonstrate how the equivalent representation of $\tilde I$ shown in Remark~\ref{thm:IW_PV_CVX} can be utilized to convexify the objective of~\eqref{prob:NWA_planning}. 

\begin{thm} \label{thm:NWAPP_cvx}

Problem~\eqref{prob:NWA_planning} can be equivalently stated as the problem 
\begin{subequations} 
\begin{align} 
\min &  \!\left\{\! \sum_{i \in \mathcal{N}} \!\left[C_i^\mathrm{O}(x_i) + I_i^\mathrm{NW}(\phi_i) \right]  
 +C^\mathrm{D}(\boldsymbol{l}^\mathrm{p})+ \frac{I}{(1 +\rho)^{\delta}} \!\right\} \\
 \mathrm{s.t.} \;& \phi_i \in \Phi_i \; \forall \; i \in \mathcal{N}\\
 &x_i \in \mathcal{X}_i (\phi_i) \; \forall \; i \in \mathcal{N} \\ 
& l^\mathrm{b}_{a,t} + \sum_{i \in \mathcal{N}} l_{a,t}^i(x_i)  \le l_a^\mathrm{p} \;\forall \; a \in \mathcal{A},\; t\in\mathcal{T} \label{eq:CVX_problem_c3}\\ 
& l_{a}^\mathrm{p} \le \overline{l}\; \forall \; a<\delta \label{eq:CVX_problem_c4} \\
& 0 \le \delta \le |\mathcal{A}| \label{eq:CVX_problem_c4.1} \\  
& \boldsymbol{l}^\mathrm{p}=\{l_a^\mathrm{p}\}_{a\in\mathcal{A}}. \label{eq:CVX_problem_c5}
\end{align} \label{eq:CVX_problem}
\end{subequations}
\end{thm}
  The proof of Thm.~\ref{thm:NWAPP_cvx} can be found in the Appendix. Since we assume that $C_i^\mathrm{O}(x_i)$, $I_i^\mathrm{NW}(\phi_i)$, and $C^\mathrm{D}(\boldsymbol{l}^\mathrm{p})$ are all convex, and $\frac{I}{(1 +\rho)^{\delta}}$ is also convex, the objective of~\eqref{eq:CVX_problem} is convex. Constraint~\eqref{eq:CVX_problem_c4}, however, introduces non-convexities to the feasible solution space of~\eqref{eq:CVX_problem}. In the next Section we show that we can decompose Problem~\eqref{eq:CVX_problem} using the DWDA and confine the non-convex constraints to be in the master problem. Then, we present an algorithm to solve the non-convex master problem by sequentially solving a small number of small linear programs.

\section{Solution technique: Dantzig-Wolfe Decomposition and Convexification} \label{sec:DWDA}
The non-convexity and dimensionality high-dimensionality of Problem~\eqref{eq:CVX_problem} may present computational challenges when using commercial solvers. As mentioned in Remark~\ref{thm:IW_PV_CVX}, the solution space is non-convex in $\bd l^\mathrm{p}$.  To illustrate the dimensionality of the problem, consider that for a time step length of 1 hour and a planning horizon of 20 years, the dimensionality of the sets $\mathcal{X}_i(\phi_i)$ ranges from roughly $175,000$ for the simplest cases (e.g., solar PV or EE) to more than half a million for the more complicated ES case. When considering all four NWAs, Problem~\eqref{eq:CVX_problem} becomes more than $1,000,000$-dimensional. 

We decompose Problem~\eqref{eq:CVX_problem} into $|\mathcal{N}|+1$ smaller subproblems to handle the dimensionality issue. Each NWA falls into a single linear problem while the demand charge and the present cost of capacity expansion are handled by the low-dimensional substation planning subproblem. Constraint~\eqref{eq:CVX_problem_c3}, however, couples all subproblems together and prevents us from independently solving them. To handle the coupling constraint, we implement the Danzig-Wolfe Decomposition Algorithm.

\subsection{Dantzig-Wolfe Decomposition}
Explicitly, the NWA subproblems are given by 
\begin{equation*}
\min_{\substack{\phi_i \in \Phi_i \\
x_i \in \bd{\mathcal{X}}_i (\phi_i)}}  \;C_i^\mathrm{O}(x_i) + I_i^\mathrm{NW}(\phi_i) +\underbrace{ \sum_{a\in\mathcal{A}} \sum_{t\in\mathcal T} \pi^1_{a,t}\cdot l^i_{a,t}(x_i)}_{\text{penalty term}}
\end{equation*}
for all $i\in\mathcal{N}$. The objective of subproblem $i$ is composed of the operation cost $C_i^\mathrm{O}(x_i)$, investment cost $I_i^\mathrm{NW}(\phi_i)$, and a term that penalizes the load $l^i_{a,t}(x_i)$ by the coefficient  $\pi^1_{a,t}$. The coefficients $ \pi^1_{a,t}$ are obtained from the dual variables of the coupling constraints~\eqref{eq:master_problem_coupling} of the master problem.  The operating and investment decisions are constrained by their respective feasible solution regions. Since $\Phi_i$ and $\bd{\mathcal{X}}_i (\phi_i)$ are convex sets and the objective of each subproblem is convex, each subproblem is a convex optimization problem. 

The master problem is given by
\begin{subequations}
\begin{align}
\min_{\delta, \boldsymbol{l}^\mathrm{p}, \lambda_k} & \left\{ \sum_{k=1}^K \lambda_k\cdot C^\mathrm{prop}_{(k)}+C^\mathrm{D}(\boldsymbol{l}^\mathrm{p})\!+\! \frac{I}{(1+\rho)^{\delta}} \right\}\\ 
\mbox{s.t. } &  l^\mathrm{b}_{a,t} + \sum_{k=1}^K\lambda_{k} \cdot  l_{a,t,(k)}^{\mathrm{prop}}  \le l_a^\mathrm{p}\;\;\; (\pi^1_{a,t})  \;\forall \; a \in \mathcal{A},\; t\in\mathcal{T} \label{eq:master_problem_coupling}\\ 
& l_{a}^\mathrm{p} \le \overline{l}\; \forall \; a<\delta \label{eq:master_problem_limit}\\
& 0 \le \delta \le |\mathcal{A}| \label{eq:master_problem_delta_lim}\\
& \boldsymbol{l}^\mathrm{p}=\{l_a^\mathrm{p}\}_{a\in\mathcal{A}} \label{eq:master_problem_peak_vec} \\
& \sum_{k=1}^K\lambda_{k} = 1 \;\;\; (\pi^2) \label{eq:master_problem_lambda_cvx} \\
& \lambda_k \ge 0 \; \forall\; k=1,\dots,K \label{eq:master_problem_lambda_pos} \\
\end{align} \label{eq:master_problem}
\end{subequations}
and its objective is to minimize a convex combination of $K$ cost proposals, peak demand charges, and the present cost of capacity expansion. The $k^\mathrm{th}$ cost proposal is defined as 
\begin{equation*}
C^\mathrm{prop}_{(k)}= \sum_{i \in \mathcal{N}} C_i^\mathrm{O}(x_{i,(k)}) + I_i^\mathrm{NW}(\phi_{i,(k)})
\end{equation*}
where $x_{i,(k)}$ and $\phi_{i,(k)}$ represent optimal operating and investment decisions for the $k^\mathrm{th}$ time the subproblem has been solved. The positive variables $\lambda_k \in [0,1]$ represent the weights assigned to each cost proposals. Constraint~\eqref{eq:master_problem_coupling} represents the coupling constraints and the load proposals are defined as
\begin{equation*}
l_{a,t,(k)}^{\mathrm{prop}}   = \sum_{i\in\mathcal{N}} l^i_{a,t}(x_{i,(k)}).
\end{equation*}
Constraints~\eqref{eq:master_problem_limit},~\eqref{eq:master_problem_delta_lim}, and~\eqref{eq:master_problem_peak_vec}  originate from Constraints~\eqref{eq:CVX_problem_c4},~\eqref{eq:CVX_problem_c4.1}, and~\eqref{eq:CVX_problem_c5}, respectively. Finally, Constraints~\eqref{eq:master_problem_lambda_cvx}  and~\eqref{eq:master_problem_lambda_pos} ensure that the sum of all $\lambda_k$'s
is equals to one (convexity constraint) and that they are all non-negative. We skip the detailed description of the well-known Danzig-Wolfe Decomposition. The interested reader is referred to~\cite{conejo2006decomposition} for an in-depth description and an implementation of the algorithm.   

\subsection{Solving the master problem}
When decomposing Problem~\eqref{eq:CVX_problem} using the DWDA, its non-convexities caused by Constraint~\eqref{eq:CVX_problem_c4} become encapsulated in the Constraint~\eqref{eq:master_problem_peak_vec} master problem. In this Section we present an algorithm to solve the master problem by solving at most $|\mathrm{A}|+1$ linear problems. Let
\begin{equation*}
\mathrm{MP}(j) 
\end{equation*}
 represent a function that fixes the variable $\delta$ to $j$ in Problem~\eqref{eq:master_problem} and solves for the variables $\bd l^\mathrm{p}$ and $\lambda_k$. A concrete interpretation of $\mathrm{MP}(j)$, is that capacity expansion happens at year $j$ and thus the peak load limit $\overline l$ is only enforced from year $1$ through $j$. Note that the function $\mathrm{MP}$ involves solving a small-scale linear problem. The number of variables in $\mathrm{MP}$ is $K+|\mathcal{A}|$ where $K$ is at most typically in the order of a few hundreds and $|\mathcal{A}|$ is, for most utility planning practices, no more than 20. 

The algorithm (shown in Algorithm~\ref{alg:MP_solve}) to solve the master problem is as follows. We sequentially solve $\mathrm{MP}(j)$ for starting with $j=0$ and increasing $j$ by one after each iteration.  If we find that $\mathrm{MP}(j)$ is greater than $\mathrm{MP}(j-1)$, we know that $\mathrm{MP}(j)$ is the optimal solution to the master problem (since its objective is convex). However, if we find that $\mathrm{MP}(j)$ is infeasible, i.e., capacity expansion cannot be delayed further than year $j-1$, we know that $\mathrm{MP}(j-1)$ is the optimal solution. Algorithm~\ref{alg:MP_solve} shows a precise description of the algorithm. 

\begin{algorithm}
	
	\KwIn{ $\mbox{MP}$, $|\bd{\mathcal A}|$ }
	\KwOut{  $\mbox{master objective value}$ }
	$j\leftarrow0$ \\ 
	\While{$j\le|\bd{\mathcal A}|$}{
		\If{$\mbox{MP}(j)$ is feasible}{
			
			\If{$j>0$ and $\mbox{MP}(j)<\mbox{MP}(j-1)$}{

				$\mbox{master objective value}\leftarrow\mbox{MP}(j-1)$  \\
				
				$j\leftarrow|\bd{\mathcal A}|+1$

			}
			
			\ElseIf{ $j =|\mathcal{A}|$ }{
				
				$\mbox{master objective value}\leftarrow \mbox{MP}(j)$  \\
				$j\leftarrow|\bd{\mathcal A}|+1$
				
			}

		}
		
		\Else{
			
			$\mbox{master objective value}\leftarrow\mbox{MP}(j-1)$  \\
			$j\leftarrow|\bd{\mathcal A}|+1$
		}
		$j\leftarrow j+1$
	}
	\caption{Solving algorithm for the master problem.}
	\label{alg:MP_solve}
	
\end{algorithm}

\section{Case Study: Non-wire alternatives for the University of Washington}
\label{sec:case_study}

The $17$ MW of load that the UW plans to add in the next 10 years could compromise the $N-1$ security of the substation and feeders that connect the campus to the Seattle City Light (SCL) distribution system. Several traditional solutions have been considered, e.g., building a new feeder in the current substation or increasing the feeder's capacity using superconducting technologies. However, these solutions are hard to implement in Seattle's dense urban environment and come at an estimated cost of over $\$60$ million. Moreover, there is an increased appetite by SCL, the Washington State government, and the UW to explore novel approaches such as NWAs. 

\subsection{Data}
Here we provide a brief description and sources of the data used in this case study.
\begin{figure} 
\centering
\includegraphics[width=.75\textwidth]{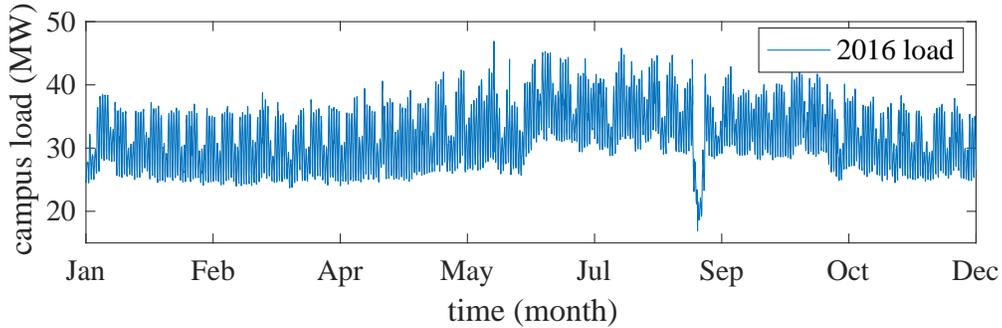}
\caption{UW Seattle campus load for the year 2016 (base load).}\label{fig:load_data}
\end{figure}
\paragraph*{Load and substation capacity}
We use one-hour electrical 2016 load data (summer peak of 48.5 MW, shown in Fig.~\ref{fig:load_data}) from the substation that serves the UW as our base load, i.e., $l_{0,t}^\mathrm{b}$, and assume that load grows at a rate of $3.5\%$ per year with respect to the base load\footnote{This rate of growth represents SCL's load growth projections of $17$ MW in 10 years.} .  The (pre-upgrade) substation capacity is $60$ MW. 

\paragraph*{Substation upgrade cost, interest rate, electricity rates, and planning horizon}
As per SCL's planning department, we assume that the cost of substation upgrades is $\$60$ million and adhere to a planning horizon of 20 years.  We assume a yearly discount rate of $7\%$ and rates based on the high-demand customer rates for the City of Seattle~\cite{SCL_tariff}.

\paragraph*{Non-wire alternatives}
Costs of EE measures are based on~\cite{brown2008us}. DR costs are based on reported values from~\cite{piette2015costs}. Costs of ES are based on substation-level lithium ion data from~\cite{lazard} and data on efficiency and degradation of ES is based on~\cite{Sarker_Optimal_2017}. Solar production profiles are based on~\cite{PFENNINGER20161251} and costs on~\cite{fu2016nrel}. Table~\ref{table:NWA_parameters} summarizes the main parameters of the NWAs considered in this study and provides their sources. In this case study, we define the operating cost of a NWA $C_i^\mathrm{O}(x_i)$ as the energy cost\footnote{Note that for net-load reduction technologies such as PV or EE, the operating cost is negative.} of the associated NWA load $l_{a,t}^i(x_i)$.

\begin{table}[]
\centering
\caption{Non-wire alternatives parameters}
\label{table:NWA_parameters}
\begin{tabular}{@{}|l|l|l|@{}}
\toprule
\rowcolor[HTML]{C0C0C0} 
\textbf{Parameter}                      & \textbf{Value}            & \textbf{Source}           \\ \midrule
\rowcolor[HTML]{EFEFEF} 
\multicolumn{3}{|l|}{\cellcolor[HTML]{EFEFEF}{\color[HTML]{000000} \textbf{Energy Efficiency}}} \\ \midrule
Investment cost function                & N/A                       & \cite{brown2008us}\footnote{Adjusted to 2017 dollars and according to the assumption that buildings in Seattle are more energy efficient than the national average (due to Seattle having some of the strictest energy codes in the nation).}                \\ \midrule
\rowcolor[HTML]{EFEFEF} 
\multicolumn{3}{|l|}{\cellcolor[HTML]{EFEFEF}\textbf{Demand response}}                          \\ \midrule
DR investment cost                      & $\$200$/kW                  &      \cite{piette2015costs}                     \\ \midrule
DR efficiency coefficient               & $1.2$                      & modeling assumption               \\ \midrule
\rowcolor[HTML]{EFEFEF} 
\multicolumn{3}{|l|}{\cellcolor[HTML]{EFEFEF}\textbf{Energy Storage}}                           \\ \midrule
Investment cost                         & $\$350$/kWh                 & \cite{lazard}                    \\ \midrule
Charge/discharge  efficiency            & $0.97/0.95$                 & \cite{Sarker_Optimal_2017}                 \\ \midrule
Degradation coefficient                 & $0.277$ kWh/kW              & \cite{Sarker_Optimal_2017}                  \\ \midrule
Energy-to-power ratio                   & $4$                         & \cite{lazard}                     \\ \midrule
\rowcolor[HTML]{EFEFEF} 
\multicolumn{3}{|l|}{\cellcolor[HTML]{EFEFEF}\textbf{Solar PV}}                                 \\ \midrule
Investment Cost                         & $\$2$/W                     & \cite{fu2016nrel}                  \\ \midrule
Production profile                      & N/A                       & \cite{PFENNINGER20161251}\footnote{Data for Seattle from the site \url{http://www.renewables.ninja}.}                    \\ \bottomrule
\end{tabular}
\end{table}

\subsection{Results}

\begin{figure} 
\centering
\includegraphics[width=.75\textwidth]{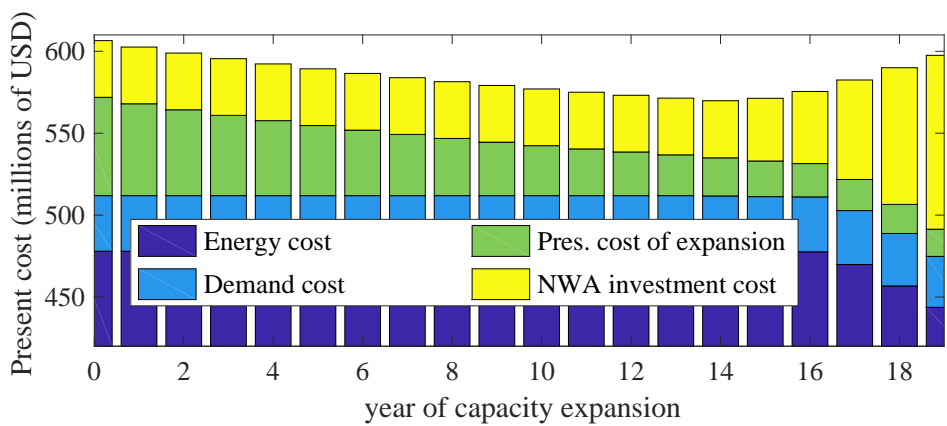}
\caption{Present cost (objective) as a function of time of expansion. The least cost solution suggest that the substation should be expanded before year 14 and that a mix of NWAs is economically viable. See Fig.~\ref{fig:3} for details on the optimal NWA mix.   }\label{fig:2}
\end{figure}

\begin{figure} 
\centering
\includegraphics[width=.75\textwidth]{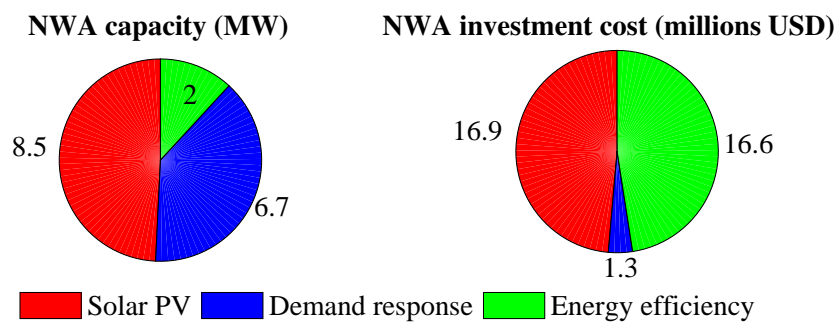}
\caption{Optimal installed NWA capacity and investment costs. Note that at \$$350$/kWh, ES is not economical. }\label{fig:3}
\end{figure}

As illustrated in Fig.~\ref{fig:2} the minimum cost is achieved when the substation expansion is delayed until year $14$ and investments include a mix of NWAs. As shown in Fig.~\ref{fig:3}, the optimal mix of NWAs include PV, DR, and EE. In this case, lithium-ion ES excluded in favor of the more economically attractive alternatives.
\begin{figure} 
\centering
\includegraphics[width=.75\textwidth]{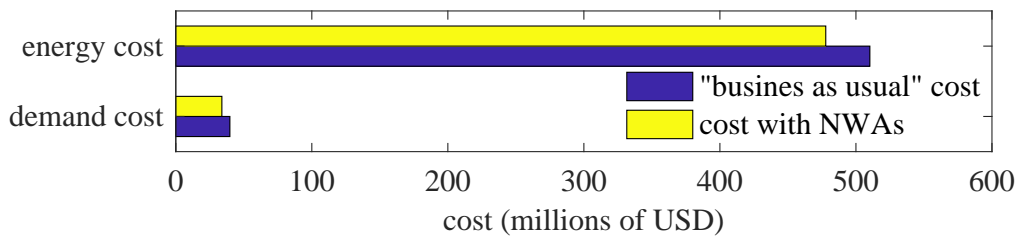}
\caption{Energy and demand costs over $20$ years with and without NWAs. NWAs are able to reduce energy costs by about $6$\% ($\$32$ million) and demand charges by about $15$\% ($\$5.8$ million). }\label{fig:E_cost_compare}
\end{figure}

It is apparent a benefit of deferring investments is that that the present cost of extension diminishes with time. However, other benefits include reduction a $6$ and $15$ percent reduction of energy and demand costs, respectively (see Fig.~\ref{fig:E_cost_compare}).

\begin{figure} 
\centering
\includegraphics[width=.75\textwidth]{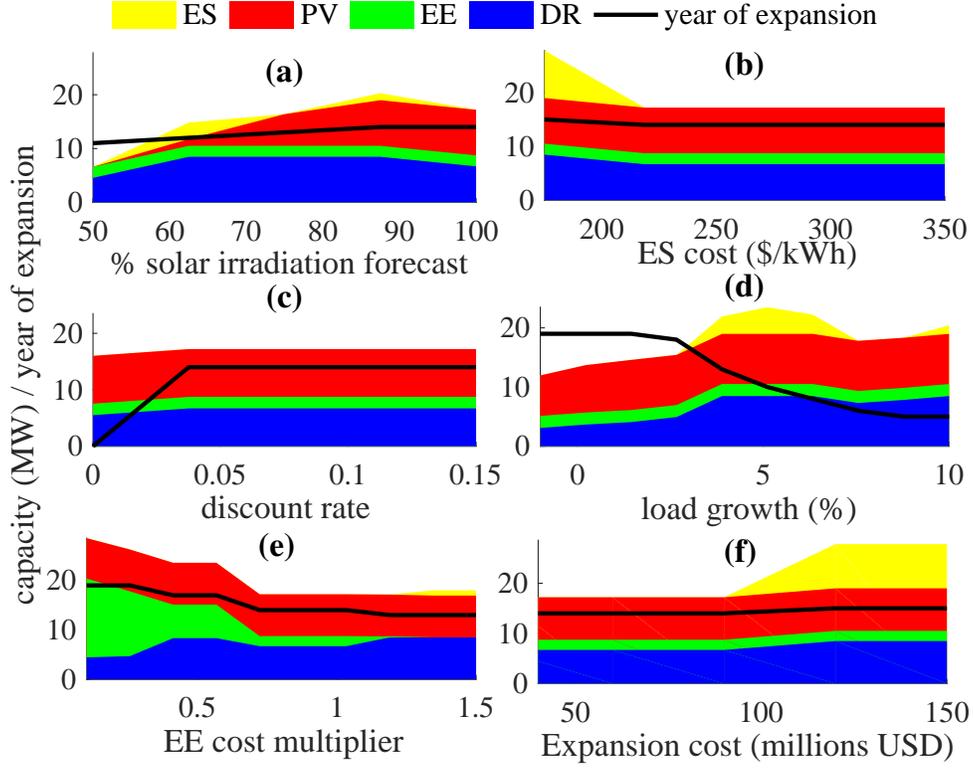}
\caption{Sensitivity of NWA installed capacities and year of expansion to four parameters: solar irradiation forecast, ES cost, the discount rate, and yearly load growth.  Plot (e) shows results when the EE cost function (i.e., the parameters $C_b^\mathrm{EE}$) scald by the EE cost multiplier. }\label{fig:sensitivity}
\end{figure}

Future extensions to this work should consider uncertainty of relevant parameters: costs, load growth projections, solar production forecasts, among others. However, considering all possible uncertainties may needlessly increase the size (and computational burden) of the model.  To find the parameters whose uncertainty may have a significant impact on the solution, we solved the NWA problem using different values for some parameters.  Fig.~\ref{fig:sensitivity} shows how the installed capacities of the four NWAs and the optimal time of expansion vary with the the value of five of the most prominent and potentially uncertain parameters. It appears that the solar irradiation forecast, load growth, and EE cost multiplier have significant impacts on the optimal NWA mix and time of expansion. Another important insight provided by Fig.~\ref{fig:sensitivity} is that ES becomes viable at less than \$$200$/kWh or when the cost of expansion is large enough.

\section{Conclusion}
\label{sec:conclusion}
We present a planning problem that determines investment and operation of distributed energy resources (DERs) and timing of capacity expansion. Considering the timing of capacity expansion has two interesting implications. First, it allows DERs to manage load and effectively act as non-wire alternatives (NWAs) to captial-intensive capacity expansion projects. Second, it makes investments in DERs more attractive by explicitly accounting for the benefit of \emph{deferring} capacity expansion investment. We formulate this problem as a large-scale non-convex optimization problem. We deal with the size of the problem by decomposing it using the Danzig-Wolfe Decomposition Algorithm. We deal with its non-convexity by further decomposing the master problem (where the non-convexities lie) into a small number of linear programs.  

Additionally, we present a realistic case study where solar photovoltaic (PV) generation, energy efficiency (EE), energy storage (ES), and demand response (DR) are considered as alternatives to substation/feeder upgrades at the University of Washington. We show that upgrades to the substation can be delayed by approximately five years by implementing PV, EE, and DR projects. At current costs, ES is not viable according to our results. Finally, we provide sensitivity analyses that suggest that uncertainty in load growth, PV generation and others should be considered in future studies. 

\appendices
\section{Proof of Theorem~\ref{thm:NWAPP_cvx}}
\label{sec:appendix}
Using the definition of $\tilde{I} $ from~\eqref{eq:IW_PV_CVX}, Problem~\eqref{prob:NWA_planning} can be written as 
\begin{equation}
\min_{\substack{\phi_i \in \Phi_i \\ x_i \in \mathcal{X}_i (\phi_i) }}  \biggl\{ \sum_{i \in \mathcal{N}} \left[C_i^\mathrm{O}(x_i) + I_i^\mathrm{NW}(\phi_i) \right]  +C^\mathrm{D}(\boldsymbol{l}^\mathrm{p}(\boldsymbol{x}))  + \min_{\delta\in\Delta } \frac{I^\mathrm{W}}{(1+\rho)^\delta}\biggr\}\label{eq:thm3_0}
\end{equation}
where $\Delta$ denotes the feasible region defined by the constraints in Problem~\eqref{eq:IW_PV_CVX}. Eq.~\eqref{eq:thm3_0} above is equivalent to 
\begin{equation}
\min_{\substack{\phi_i \in \Phi_i \\ x_i \in \mathcal{X}_i (\phi_i) }} \min_{\delta\in\Delta } \biggl\{ \sum_{i \in \mathcal{N}} \left[C_i^\mathrm{O}(x_i) + I_i^\mathrm{NW}(\phi_i) \right]   +C^\mathrm{D}(\boldsymbol{l}^\mathrm{p}(\boldsymbol{x}))  +  \frac{I}{(1+\rho)^\delta}\biggr\} \label{eq:thm3_1}.
\end{equation}
Problem~\eqref{eq:thm3_1} is a nested optimization problem that whose inner variable is $\delta$ and its outer variables are $x_i$ and $\phi_i$. As shown in reference~\cite{boyd2004convex} all variables can be minimized simultaneously in an equivalent problem. Thus~\eqref{prob:NWA_planning} is equivalent to~\eqref{eq:CVX_problem} where Constraints~\eqref{eq:CVX_problem_c3} are employed to implement the functions $l^\mathrm{p}_a(\boldsymbol{x}) = \max_{t\in\mathcal{T}}\{l_{a,t}(\boldsymbol{x})\}$ via half-planes. 
\QEDB

\bibliographystyle{IEEEtran}
\bibliography{./sample.bib}

\end{document}